\documentclass[letterpaper, 10 pt, journal, twoside]{IEEEtran}

\usepackage[utf8]{inputenc}
\usepackage{amsmath}
\usepackage{amssymb}
\usepackage{enumerate}
\usepackage{cite}
\usepackage[allowmove]{url}
\usepackage{color}
\usepackage{graphicx}          

\usepackage[ruled,vlined]{algorithm2e}

\usepackage{mathtools}

\def\R{{\mathbb R}}
\def\N{{\mathbb N}}

\def\sign{\,\mathrm{sign}}

\def\C{\mathcal{C}}
\def\F{\mathcal{F}}
\def\M{\mathcal{M}}

\def\W{\mathcal{W}}
\def\X{\mathcal{X}}
\def\Oo{\mathbf{0}}
\def\one{\mathbf{1}}
\def\b{\mathbf{b}}

\def\x{\mathbf{x}}

\def\m{\mathbf{m}}

\def\id{\mathbf{I}}

\def\fep{{\lfloor\varepsilon\rfloor}}
\def\ep{\varepsilon}
\def\f{\mathfrak{f}}

\newtheorem{rem}{Remark}
\newtheorem{lem}{Lemma}
\newtheorem{prob}{Problem}

\newtheorem{thm}{Theorem}

\title{
Differentiator for Noisy Sampled Signals with\\ Best Worst-Case Accuracy* 
}

\author{Hernan Haimovich, Richard Seeber, Rodrigo Aldana-L\'opez, and David G\'omez-Guti\'errez
\thanks{*Work supported by FONCYT, Argentina, under grant PICT 2018-1385, and by Christian Doppler Research Association, Austrian Federal Ministry for Digital and Economic Affairs, and National Foundation for Research, Technology and
Development.}
\thanks{Hernan Haimovich is with Centro Internacional Franco-Argentino de 
  Ciencias de la Informaci\'on y de Sistemas (CIFASIS)
  CONICET-UNR, 2000 Rosario, Argentina.
        (e-mail: haimovich@cifasis-conicet.gov.ar)}%
\thanks{Richard Seeber is with the Christian Doppler Laboratory for Model Based Control of Complex Test Bed Systems, Institute of Automation and Control, Graz University of Technology, Graz, Austria.
        (e-mail: richard.seeber@tugraz.at)}%

\thanks{Rodrigo Aldana-L\'opez is with the Department of Computer Science and Systems Engineering, University of Zaragoza, Zaragoza, Spain.  
        (e-mail: rodrigo.aldana.lopez@gmail.com)}%
        
\thanks{
David G\'omez-Guti\'errez is with Intel Labs, Intel Corporation and with Tecnologico de Monterrey, Jalisco, Mexico. 
        (e-mail: david.gomez.g@ieee.org)}%
\thanks{\textcolor{red}{This is accepted version of the manuscript: Hernan Haimovich, Richard Seeber, Rodrigo Aldana-L\'opez, and David G\'omez-Guti\'errez, ``Differentiator for Noisy Sampled Signals with Best Worst-Case Accuracy". IEEE Control Systems Letters. DOI: 10.1109/LCSYS.2021.3087542. 
\textbf{Please cite the publisher's version}. For the publisher's version and full citation details see:
\url{https://doi.org/10.1109/LCSYS.2021.3087542}. 
}}
\thanks{``© 2021 IEEE.  Personal use of this material is permitted.  Permission from IEEE must be obtained for all other uses, in any current or future media, including reprinting/republishing this material for advertising or promotional purposes, creating new collective works, for resale or redistribution to servers or lists, or reuse of any copyrighted component of this work in other works.”}
}
\begin{document}

\maketitle
\thispagestyle{empty}
\pagestyle{empty}

\begin{abstract}
This paper proposes a differentiator for sampled signals with bounded noise and bounded second derivative.
It is based on a linear program derived from the available sample information and requires no further tuning beyond the noise and derivative bounds.
A tight bound on the worst-case accuracy, i.e., the worst-case differentiation error, is derived, which is the best among all causal differentiators and is moreover shown to be obtained after a fixed number of sampling steps.
Comparisons with the accuracy of existing high-gain and sliding-mode differentiators illustrate the obtained results.
\end{abstract}

\section{Introduction}
\label{sec:intro}

\IEEEPARstart{E}{stimating} in real-time the derivatives of a signal affected by noise is a fundamental problem in control theory and continues to be an active area of research, see, e.g., the special issue~\cite{Reichhartinger2018SpecialDifferentiators}, the comparative analysis~\cite{RasoolMojallizadeh2021Discrete-timeAnalysis}, and the references therein.
Differentiators are often used, for instance, for state estimation~\cite{Shtessel2014ObservationObservers}, Proportional-Derivative controllers, fault detection~\cite{Rios2015fault,Efimov2012ApplicationDetection}, and unknown input observers~\cite{Bejarano2010HighInputs}.

Popular existing methods for differentiation include linear high-gain observers~\cite{Vasiljevic2008ErrorObservers}, linear algebraic differentiators~\cite{Mboup2007AControl,Othmane2020AnalysisDisturbances}, and sliding mode differentiators~\cite{Levant1998RobustTechnique,Levant2003}.
These differ in terms of their convergence properties; while high-gain differentiators converge exponentially~\cite{Vasiljevic2008ErrorObservers}, algebraic and sliding-mode differentiators exhibit convergence in finite or fixed time~\cite{Mboup2007AControl,Levant1998RobustTechnique}, and they converge exactly for different classes of noise-free signals~\cite{Levant2003,Othmane2020AnalysisDisturbances}.
With measurement noise, the accuracy, i.e., the achievable worst-case differentiation error, is limited for all differentiators.
Linear differentiators may be tuned to minimize the bound on the differentiation error when the noise amplitude and a bound on the derivative's Lipschitz constant are known~\cite{Vasiljevic2008ErrorObservers}, whereas the tuning of sliding mode differentiators only requires knowledge about the Lipschitz constant but not about the noise~\cite{Levant1998RobustTechnique,Fraguela2012}.

In practice, differentiation is typically performed on a digital computer using sampled signals. Hence, the use of continuous-time differentiators requires discretization, which is particularly challenging for sliding-mode differentiators because an explicit (forward) Euler discretization may lead to reduced accuracy, numerical chattering, and even instability \cite{Polyakov2019ConsistentSystems,Levant2013OnControl}.
Several techniques for that purpose have therefore been proposed, cf. \cite{RasoolMojallizadeh2021Discrete-timeAnalysis,Carvajal-Rubio2021ImplicitDifferentiators,Koch2018DiscreteDifferentiators}.
In any case, the inherent performance limitations of continuous-time differentiators cannot be surpassed in the discrete domain via discretization.

The present paper proposes a differentiator that considers the information available in the sampled signal in the form of a linear program.
This approach also yields upper and lower bounds for the derivative, similar to interval observers~\cite{Mazenc2011IntervalDisturbances}.
Interval observers, however, have seldom been applied to differentiation, see e.g., ~\cite{Guerra2017IntervalAccuracy}, and they are limited in terms of accuracy by their underlying observer.

In contrast to other observers, the present approach is shown to have the best possible worst-case accuracy among all causal differentiators.
This best possible worst-case accuracy is shown to be achieved using a fixed number of samples, thus providing a limit on the computational complexity of the linear program and guaranteeing convergence in a fixed time, similarly to algebraic and some sliding-mode differentiators.
Moreover, implementing the algorithm only requires knowledge of the derivative's Lipschitz constant and the noise bound but, unlike other differentiators, yields such an estimate without requiring any further tuning.

\textbf{Notation:} $\N$, $\N_0$, $\R$ and $\R_{\ge 0}$ denote the positive and nonnegative integers, and the reals and nonnegative reals, respectively. If $\alpha\in\R$, then $|\alpha|$ denotes its absolute value. 
For $x,y\in\R^n$, inequalities and absolute value operate componentwise, so that $|x|\in\R_{\ge 0}^n$ denotes the vector with components $|x_i|$, and $x\le y$ the set of inequalities $x_i \le y_i$, for $i=1,\ldots,n$. For a (differentiable) function $f : D\subset\R \to \R$, $f^{(i)}$ denotes its $i$-th order derivative. For $a\in\R$, the greatest integer not greater than $a$ is denoted by $\lfloor a \rfloor$. The symbols $\Oo$, $\id$ and $\one$ denote the zero vector or matrix, the identity matrix and a vector all of whose components are equal to one, respectively.

\section{Problem Statement and Formulation}
\label{sec:preliminaries}

\subsection{Problem statement}

Consider a differentiable function $\f : \R_{\ge 0} \to \R$ of which we know that its derivative $\f^{(1)}$ is globally Lipschitz continuous with Lipschitz constant $L$, i.e.
\begin{align}
  \label{eq:fddotbnd}
  |\f^{(2)}(t)| &\le L,\quad \text{for almost all }t\in\R_{\ge 0}.
\end{align}
Suppose that a noisy measurement $m_k$ of the value $\f(kT)$ becomes available at each time instant $t_k = kT$, with $k \in \N_0$, and that a bound $N$ on the noise is known, so that
\begin{align}
  \label{eq:basemnf}
  m_j &= \f_j + \eta_j, & \f_j &:= \f(jT), &|\eta_j| &\le N,
\end{align}
for $j = 0,1,\ldots,k$.

The problem to be addressed is to design an algorithm that, at every time instant $t_k$ when a new measurement $m_k$ becomes known, gives all available information on the current value 
\begin{align}
  \label{eq:dersample}
  \f_k^1 := \f^{(1)}(t_k),\quad t_k := kT,
\end{align}
of the derivative of $\f$. The more specific problem to be solved is as follows.
\begin{prob}
  \label{prob:main}
  Devise an algorithm that, given the constants $L\ge 0$ (bound on 2nd-order derivative), $N\ge 0$ (noise bound), $T>0$ (sampling period) and $k\in\N$, provides all possible values of $\f^1_k$ based on knowledge of the bound \eqref{eq:fddotbnd} and the measurements $m_j$ for $j=0,1,\ldots,k$. 
\end{prob}

\subsection{Possible values for the derivative}

Let $\F_k^1(\m_k)$ denote the set of possible values for $\f_k^1 = \f^{(1)}(t_k)$ that are consistent with the bound \eqref{eq:fddotbnd} and the measurements 
\begin{align}
    \label{eq:defmk}
    \m_k &:= [m_0,\ m_1,\ \ldots,\ m_k]^T
\end{align}
that satisfy~\eqref{eq:basemnf}. The set $\F_k^1(\m_k)$ can be defined as
\begin{align*}
    \F^1_k(\m_k) := \{ \f_k^1 \in \R : \exists \f(\cdot) \text{ satisfying }\eqref{eq:fddotbnd}, \eqref{eq:basemnf}, \eqref{eq:dersample}\}.
\end{align*}
The set $\F_k^1(\m_k)$ is convex and hence, whenever nonempty, it will have the form of an interval. Problem~\ref{prob:main} can thus be posed as finding the extreme values
\begin{align}
    \label{eq:trueworst}
    \overline{\F}_k^1(\m_k) &:= \sup \F_k^1(\m_k), &\underline{\F}_k^1(\m_k) &:= \inf \F_k^1(\m_k).
\end{align}
For future reference, define $\X_k(\f)$ as the vector
\begin{align}
    \label{eq:defXkf}
    \X_k(\f) := [\f(t_0),\ldots,\f(t_k), \f^{(1)}(t_0),\ldots,\f^{(1)}(t_k)]^T.
\end{align}

\subsection{Samples and measurements}

Since the derivative $\f^{(1)}$ is globally Lipschitz continuous, then $\f^{(2)}$ exists almost everywhere and
\begin{align*}
  \f^{(1)}(\tau) &= \f^{1}_{j} - \int\limits_{\tau}^{t_j} \f^{(2)}(s)  ds,  
  & \f(t) &= \f_j - \int\limits_{t}^{t_j} \f^{(1)}(\tau)  d\tau.
\end{align*}

From these expressions one can obtain the bounds
\begin{align}
  \label{eq:ff1const}
  |\f_{j-1} - \f_j + \f_j^1 T| &\le L \frac{T^2}{2},\\
  \label{eq:f1const}
  |\f_j^1 - \f_{j-1}^1| &\le LT,\\
\text{and from~\eqref{eq:basemnf}, also \quad}
  \label{eq:fmbnd}
  |\f_j - m_j| &\le N.
\end{align}
At time $t_k = kT$, every function $\f$ that satisfies the bound~\eqref{eq:fddotbnd} for almost all $t \in [0,t_k]$ will have associated values $\X_k(\f)$ that must satisfy the constraints~\eqref{eq:ff1const}--\eqref{eq:f1const} for $j=1,2,\ldots,k$. In addition, given the noise bound $N$, the corresponding measurements must satisfy~\eqref{eq:fmbnd} for $j=0,1,\ldots,k$. 

\section{Main Results}
\label{sec:mainresults}

\subsection{Derivation of the proposed differentiator}

Consider a vector $\x_k$ of $2k+2$ \emph{optimization variables}
\begin{align}
  \label{eq:defxk}
  \x_k &:= [(f_{0:k})^T,\ (f_{0:k}^1)^T]^T \in \R^{2k+2}\\
  f_{0:k} &:= [f_0,\ f_1,\ \ldots,\ f_k]^T, \quad
  f_{0:k}^1 := [f_0^1,\ f_1^1,\ \ldots,\ f_k^1]^T,\notag
\end{align}
where $f_i$ and $f^1_i$ model possible (hypothetical) values for $\f_i$ and $\f^1_i$, respectively.
For every $\m\in\R^{k+1}$, consider the set
\begin{align}
  \label{eq:cset}
  &\C_k(\m) := \{ \x \in \R^{2k+2} : |A_k \x + M_k \m| \le \b_k \},\\
  &A_k = 
  \begin{bmatrix}
    \Oo & D_k\\
    D_k & -T[\Oo\ \id]\\
    \id & \Oo
  \end{bmatrix},\quad
  M_k = 
  \begin{bmatrix}
    \Oo \\ \Oo \\ \id
  \end{bmatrix},\quad
  \b_k = 
  \begin{bmatrix}
    LT \one \\ L\dfrac{T^2}{2} \one \\ N \one
  \end{bmatrix}
  \notag
\end{align}
where $D_{k} \in \R^{k\times k+1}$ is a Toeplitz matrix with first row and column given by $[-1, 1, 0, \ldots, 0]$ and $[-1, 0, \ldots, 0]^T$, respectively.
The set $\C_k(\m)$ is defined so that whenever a function $\f$ satisfies~\eqref{eq:fddotbnd} and produces the measurements $\m_k$ satisfying \eqref{eq:basemnf}, then $\X_k(\f) \in \C_k(\m_k)$. This is so because 
the rows of the matrices $A_k,M_k$ and the vector $\b_k$ are grouped into 3 blocks, of $k$, $k$ and $k+1$ rows, where the first block corresponds to \eqref{eq:f1const}, the second to \eqref{eq:ff1const}, and the third to \eqref{eq:fmbnd}.
\begin{rem}
\label{rem:bgg}
Although $\X_k(\f) \in \C_k(\m_k)$ holds for all admissible functions $\f$ and corresponding measurements $\m_k$, given an arbitrary vector $\m_k\in \R^{k+1}$ with nonempty $\C_k(\m_k)$ it may not be true that a function $\f$ exists satisfying \eqref{eq:fddotbnd}--\eqref{eq:basemnf} with $\X_k(\f) \in \C_k(\m_k)$ (see the Appendix for a counterexample).
\end{rem}
The proposed differentiator provides an estimate $\hat{\f}_k^1$ for the derivative $\f_k^1 = \f^{(1)}(t_k)$ by solving the optimization problems \eqref{eq:bndmax}--\eqref{eq:bndmin} and computing \eqref{eq:estimate}:
\begin{subequations}
  \label{eq:bounds}
  \begin{align}
    \label{eq:bndmax}
    \overline{f}_k^1(\m_k) &:= \max \{c_k^T \x_k \ : \x_k \in \C_k(\m_k)\},\\
    \label{eq:bndmin}
    \underline{f}_k^1(\m_k) &:= \min \{c_k^T \x_k \ : \x_k \in \C_k(\m_k)\},\\
    \label{eq:estimate}
    \hat\f_k^1 &:= \left( \overline{f}_k^1(\m_k) + \underline{f}_k^1(\m_k) \right)/2,
  \end{align}
\end{subequations}
with
$
  c_k = [0,\ \ldots,\ 0,\ 1]^T \in \R^{2k+2}.
$
Note that $c_k^T \x_k = f_k^1$, according to \eqref{eq:defxk}.
From Remark~\ref{rem:bgg}, it follows that the set of possible values for the derivative of $\f$ at time $t_k$, namely $\F_k^1(\m_k)$, satisfies $\F_k^1(\m_k) \subseteq [\underline{f}_k^1(\m_k),\overline{f}_k^1(\m_k)]$ and thus 
\begin{align}
    \label{eq:truelpineq}
    \underline{f}_k^1(\m_k) \le \underline{\F}_k^1(\m_k) \le \overline{\F}_k^1(\m_k) \le \overline{f}_k^1(\m_k).
\end{align}
The set $\C_k(\m_k)$ is defined by linear inequalities in the optimization variables, for every $\m_k$. Thus, \eqref{eq:bndmax} and \eqref{eq:bndmin} are linear programs; the only information required to implement them are the values $L$, $N$, $T$ and the measurements $\m_k$, obtained up to $t_k$. The proposed estimate $\hat\f_k^1$ yields the smallest worst-case distance to any value within  $[\underline{f}_k^1(\m_k),\overline{f}_k^1(\m_k)]$.

The computational complexity of the linear programs increases with increasing $k$. A fixed number of samples $\hat{K}+1$ can be considered to limit the complexity as summarized in Algorithm \ref{algo:differentiator}, which is meant to be executed at every time instant.  The next section provides a way to choose $\hat{K}$ by studying the worst-case accuracy of the differentiator and showing that a finite $K$ can be computed such that for all $\hat{K}\geq K$ the same worst-case accuracy  is obtained.

\begin{algorithm}
\label{algo:differentiator}
\SetAlgoLined
\SetKwInput{Input}{input}
\SetKwInOut{Return}{return}

\Input{$L$, $N$, $T$, $\hat{K}$, $\m_k$}

Set $\underline k := \min\{k, \hat{K}\}$

Set $\m_{\underline k}^{k} := [m_{k-\underline k},\dots,m_{k-1},m_k]^T\in\mathbb{R}^{\underline k+1}$

Set $A_{\underline k}$, $M_{\underline k}$, and $\b_{\underline k}$ as in \eqref{eq:cset} using $L$, $N$, $T$.

Set $c_{\underline k}:=[0,\ \ldots,\ 0,\ 1]^T \in \R^{2\underline k+2}$

Solve $\overline{f}_k^1 := \max \left\{c_{\underline k}^T \x \ :  |A_{\underline k} \x + M_{\underline k} \m_{\underline k}^{k}| \le \b_{\underline k} \right\}$

Solve $\underline{f}_k^1 := \min \left\{c_{\underline k}^T \x \ :  |A_{\underline k} \x + M_{\underline k} \m_{\underline k}^{k}| \le \b_{\underline k} \right\}$

\Return{$\hat\f_k^1 := \left( \overline{f}_k^1 + \underline{f}_k^1 \right)/2$}

\caption{Estimation of $\f^{(1)}(kT)$, based on $\hat{K}+1$ noisy measurements, using linear programming.}
\end{algorithm}

\subsection{Differentiator convergence and worst-case accuracy}

A measure of the accuracy of the differentiator is given by the difference between the upper and lower bounds that can be ensured on the derivative
\begin{equation}
  \label{eq:width}
  w_k(\m_k) := \overline{f}_k^1(\m_k) - \underline{f}_k^1(\m_k).
\end{equation}
With the differentiator output $\hat\f_k^1$ suggested above, the differentiator error is then bounded by $w_k(\m_k)/2 \ge 0$ if $\C_k(\m_k) \neq \emptyset$.
A related quantity is the difference of actual worst-case derivative bounds
\begin{equation}
  \label{eq:truewidth}
  \W_k(\m_k) := \overline{\F}_k^1(\m_k) - \underline{\F}_k^1(\m_k),
\end{equation}
which according to \eqref{eq:trueworst} correspond to the best possible accuracy obtainable from the measurements $\m_k$.

Let $\M_k(L,N,T)$ denote the set of all possible measurements $\m_k$ that could be obtained for functions satisfying~\eqref{eq:fddotbnd} with additive measurement noise bound $N$:
\begin{align*}
  \M_k(L,N,T) := \big\{ \m_k &\in\R^{k+1} : \exists \f(\cdot), \text{ \eqref{eq:fddotbnd} and~\eqref{eq:basemnf} hold}\big\}.
\end{align*}
Consider the obtained and the best possible accuracy over all possible measurements, i.e., their worst-case values,
\begin{align}
  \label{eq:wc-diff}
  \bar{w}_k(L,N,T) &:= \textstyle\sup_{\m \in \M_k(L,N,T)} w_k(\m), \\
  \label{eq:wc-true}
  \bar{\W}_k(L,N,T) &:= \textstyle\sup_{\m \in \M_k(L,N,T)} \W_k(\m).
\end{align}
Clearly, $\bar w_k(L,N,T) \ge \bar{\W}_k(L,N,T)$ according to \eqref{eq:truelpineq}.
Also, no causal differentiator can achieve a better worst-case accuracy than $\bar{\W}_k(L,N,T)$ due to \eqref{eq:trueworst}.

Our main result is the following.
\begin{thm}
  \label{thm:zero-meas-worst}
  Given positive $L$, $N$, and $T$, the accuracies $w_k(\m)$, $\bar w_k(L,N,T)$ obtained with the differentiator~\eqref{eq:bounds} and the best possible accuracies $\W_k(\m)$, $\bar \W_k(L,N,T)$, as defined in \eqref{eq:width}, \eqref{eq:wc-diff}  and \eqref{eq:truewidth}, \eqref{eq:wc-true}, respectively, satisfy:
  \begin{enumerate}[a)]
    \item \label{item:finconv}$w_k(\Oo) = 2\bar h_k$, with
    \begin{align}
      \bar h_k &= h_o(\underline k), \ \underline k = \min\{k,K\}, \ h_o(\ell) := \frac{1}{2}LT\ell + \frac{2N}{T\ell},\notag\\
      K &:= 
      \begin{cases}
        Q &\text{if } Q^2+Q\geq \frac{4N}{LT^2},\\
        Q + 1 & \text{otherwise,}
      \end{cases} \quad
      Q := \left\lfloor \frac{2}{T} \sqrt{\frac{N}{L}} \right\rfloor.\notag
    \end{align}
    \item $w_{k}(\Oo) \ge w_{k+1}(\Oo)$ for all $k\in\N$;\label{item:noninc}    
    \item $w_k(\Oo) = \W_k(\Oo)$;\label{item:0lp0fun}
    \item $w_k(\Oo) = \bar{w}_k(L,N,T) = \bar{\W}_k(L,N,T)$.\label{item:0worst}
  \end{enumerate}
\end{thm}

Items~\ref{item:finconv}) and~\ref{item:noninc}) state that the sequence $\{w_k(\Oo)\}$ is nonincreasing and converges to a limit in $K$ samples, and give an expression for both, the number of samples and the limit value $\bar{h}_K$.
Item~\ref{item:0lp0fun}) states that, when all measurements equal zero, the accuracy obtained is identical to the true, best possible accuracy among causal differentiators formulated in Problem~\ref{prob:main}.
Item~\ref{item:0worst}) shows that the zero-measurement case is actually the worst over all possible measurements. This means that the proposed differentiator's worst-case accuracy is thus the best among all causal differentiators.
Note that these results are very powerful because $\x_k \in \C_k(\m_k)$ does not imply that $\x_k = \X_k(\f)$ for some $\f$ satisfying \eqref{eq:fddotbnd}--\eqref{eq:basemnf}, as stated in Remark~\ref{rem:bgg}. 

Since the best worst-case accuracy is achieved after a fixed number of $K$ sampling steps and then stays constant, considering more (older) measurements does not improve the worst-case performance. With this insight, Theorem~\ref{thm:zero-meas-worst} ensures that Algorithm~\ref{algo:differentiator} with $\hat{K}\geq K$ provides the best worst-case accuracy among all causal differentiators.
Particularly, if $N < L T^2/4$, then $\hat{K} = K = 1$ can be chosen; the linear programs may then be solved explicitly, yielding the differentiator $\hat\f_k^1 = (m_{k} - m_{k-1})/T$ as a special case.

\section{Proof of Theorem \ref{thm:zero-meas-worst}}
\label{sec:proof}

The proof strategy for Theorem \ref{thm:zero-meas-worst} is to first study the case corresponding to the noise bound $N=1$ and the sampling period $T=1$, and then show how the general case can be obtained from this. 
\subsection{The case $N=T=1$}
To begin, consider the case $(N,T,L)=(1,1,\tilde{L})$ with $\tilde{L}:=4/\varepsilon^2$. Using these parameters, the quantities in the statement of Theorem~\ref{thm:zero-meas-worst} become
\begin{equation}
    \begin{aligned}
    \label{eq:definitions}
 K &= \begin{cases}
        \lfloor\varepsilon\rfloor &\text{if } \lfloor\varepsilon\rfloor^2+\lfloor\varepsilon\rfloor \ge \ep^2,\\
        \lfloor\varepsilon\rfloor + 1 & \text{otherwise,}
      \end{cases}\\    
      \bar a_k &= h(\underline k), \ \underline k = \min\{k,K\}, \ h(\ell)= \frac{2\ell}{\varepsilon^2}+\frac{2}{\ell},
    \end{aligned}
\end{equation}
where we have used $\overline a_k$ and $h(\cdot)$ to denote $\bar h_k$ and $h_o(\cdot)$ corresponding to $(N,T,L)=(1,1,4/\varepsilon^2)$.

The following lemma establishes some properties of the sequence $\{\overline{a}_k\}$ that will be required next.
\begin{lem}
\label{le:sequence_a}
Consider \eqref{eq:definitions}, and let
  \begin{align*}
    a_k &:= \min_{\ell \in \{1,\ldots,k\}} h(\ell).
  \end{align*}
  Then, the following statements are true:
  \begin{enumerate}[a)]
  \item $h(\ell)$ is strictly decreasing for $\ell\leq \fep$ and strictly increasing for $\ell\geq \fep+1$.
  \item If $k\leq \fep$ then $a_{k}=h(k)$.
  \item If $k\geq \fep+1$ then $a_k = a_K$.
  \item $\overline{a}_k = a_k$ for all $k\in\N$.
  \end{enumerate}
\end{lem}
\begin{IEEEproof}
The derivative of the function $h$ is $h^{(1)}(s) = \frac{2}{\varepsilon^2} - \frac{2}{s^2}$,
so that $h^{(1)}(\varepsilon) = 0$, $h^{(1)}(s) < 0$ for $s\in(0,\varepsilon)$, and $h^{(1)}(s) > 0$ for $s\in (\varepsilon,\infty)$. Therefore, $h$ is strictly decreasing within the interval $(0,\varepsilon]$ and strictly increasing within $[\varepsilon,\infty)$. Since $\fep \le \varepsilon < \fep + 1$, then item a) is established. Item b) then follows straightforwardly from the definition of $a_k$.

For item c), note that for $k \ge \fep + 1 > \varepsilon$, from item a) we must have $a_k=\min\{h({\fep}),h({\fep+1})\}$. Consider
\begin{align*}
    h(\fep) - h(\fep+1) &=  \frac{2}{\fep} - \frac{2}{\fep+1}-\frac{2}{\ep^2}  
  \end{align*}
  If this difference is nonpositive, which happens if
  $
  \fep^2+\fep \ge \ep^2,
  $
  then $h({\fep+1}) \ge h(\fep)$ will hold. Observing \eqref{eq:definitions}, then
$a_k = \min\{h(\fep),h({\fep+1})\} = a_K$. 

Finally, d) follows by combining Lemma \ref{le:sequence_a}b) and c).
\end{IEEEproof}

Let $S_k:=\{\ell\in\N: 1\leq \ell\leq k-\underline{k}\}$. Consider a function $\tilde{\f} : [0,k]\to\R$ defined as follows
\begin{align}
  \label{eq:worse_f_1}
  \tilde{\f}(t) &:= 
    \frac{2(t-k)^2}{\ep^2}+\overline{a}_k(t-k)+1\ \text{for } t\in[k-\underline{k},k]\ \\
  \label{eq:worse_f_2}
   \tilde{\f}(t) &:=\tilde{\f}^{(1)}(\ell)(t-\ell)^2 + \tilde{\f}^{(1)}(\ell)(t-\ell) - 1
  \end{align}
for $t\in[\ell-1,\ell)$ with $\ell\in S_k$ and where
\begin{equation}
    \label{eq:f1rightcont}
    \tilde \f^{(1)}(\ell) := \lim_{t\to \ell^+} \tilde \f^{(1)}(t).
\end{equation}
It is clear that $\tilde \f$ satisfies $\lim_{t\to k^-} \tilde \f^{(1)}(t) = \overline a_k$. The following lemma establishes that $\tilde \f$ is continuously differentiable, its derivative has a global Lipschitz constant $\tilde L$, and satisfies $\X_k(\tilde \f) \in \C_k(\Oo)$.
\begin{lem}
\label{lem:prop_ftilde}
  Let $\tilde{\f}:[0,k]\to\mathbb{R}$ be defined by \eqref{eq:worse_f_1}--\eqref{eq:f1rightcont}. Then,
  \begin{enumerate}[a)]
      \item $\tilde{\f}$ is continuous, $\tilde{\f}(k)=1$, and $\tilde{\f}(\ell)=-1$ for $\ell\in S_k$.\label{item:prop_f}
      \item $\tilde{\f}^{(1)}$ is continuous in $(0,k)$, and $\tilde{\f}^{(1)}(\ell-1) = -\tilde{\f}^{(1)}(\ell)$ for every $\ell\in S_k$.\label{item:prop_df}
      \item $|\tilde{\f}^{(1)}(k-\underline{k})| \le 2/\ep^2$ for $k > \underline{k}$\label{item:prop_dfboundk}
      \item $\tilde{\f}(\ell)\in[-1,1]$ for every $\ell\in\N_0$ with $\ell \leq k$.\label{item:prop_fbound}
      \item $|\tilde{\f}^{(2)}(t)|\leq \tilde{L}=4/\varepsilon^2$ for almost every $t\in[0,k]$. \label{item:prop_dfbound}
  \end{enumerate}
\end{lem}
\begin{IEEEproof}
\ref{item:prop_f}) The fact that $\tilde{\f}(k)=1$ follows directly from \eqref{eq:worse_f_1}. Also, since $\overline{a}_k = h(\underline k)$, then
$$
\tilde{\f}(k-\underline{k})=\frac{2\underline{k}^2}{\ep^2} + \left(\frac{2\underline{k}}{\ep^2}+\frac{2}{\underline{k}}\right)(-\underline{k})+1 = -1.
$$
Note that, by definition, $\tilde \f$ is continuous in the intervals $[\ell-1,\ell)$ for all $\ell \in S_k$ and also in $[k-\underline k,k]$. From \eqref{eq:worse_f_2}, for $\ell\in S_k$ one has
$\lim_{t\to\ell^-}\tilde{\f}(t) = -1 = \tilde \f(\ell)$. Thus, $\tilde{\f}(\ell)=-1$ for any $\ell\in S_k$ and it follows that $\tilde{\f}$ is continuous in $[0,k]$.

\ref{item:prop_df}) From \eqref{eq:worse_f_2} we obtain:
\begin{equation}
\label{eq:worse_df_2}
\tilde{\f}^{(1)}(t) = 2\tilde{\f}^{(1)}(\ell)(t-\ell)+\tilde{\f}^{(1)}(\ell), \quad t\in[\ell-1,\ell)
\end{equation}
Hence, \eqref{eq:worse_df_2} gives $\lim_{t\to\ell^-} \tilde{\f}^{(1)}(t) = \tilde{\f}^{(1)}(\ell)$ which according to \eqref{eq:f1rightcont} leads to continuity of $\tilde{\f}^{(1)}(t)$ for every $t\in S_k$. Continuity of $\tilde{\f}^{(1)}$ within $(0,k]$ then follows similarly as in the proof of item~\ref{item:prop_f}. Finally, note that evaluating at $t=\ell-1$ in \eqref{eq:worse_df_2} it follows that $\tilde{\f}^{(1)}(\ell-1)=-\tilde{\f}^{(1)}(\ell)$ for every $\ell\in S_k$.

\ref{item:prop_dfboundk}) From \eqref{eq:definitions}, \eqref{eq:worse_f_1} and the definition $\overline{a}_k = h(\underline k)$, it follows that for $k > \underline{k} = K$
\begin{equation}
\label{eq:critical_df}
\tilde{\f}^{(1)}(k-\underline{k}) = \frac{2}{K}-\frac{2K}{\ep^2}.
\end{equation}
Multiplying $\tilde \f^{(1)}(k-\underline{k})$ by $K\ep^2/2>0$, the inequalities
$
-K \le \ep^2 - K^2 \le K
$
have to be proven.
Consider first the case $K = \fep$.
Then, $K \le \ep$, i.e., $\ep^2 - K^2 \ge 0$, and the upper bound remains to be proven.
Since $\fep^2+\fep \ge \ep^2$ holds in this case, one has
\begin{equation*}
    \ep^2 - \fep^2 \le \ep^2 - (\ep^2  - \fep) = \fep = K.
\end{equation*}
Consider now the case $K = \fep + 1$.
Since $K > \ep$, it suffices to show the lower bound.
It is obtained from
\begin{equation*}
    \ep^2 - (\fep +1)^2 = \ep^2 - \fep^2 - 2\fep - 1 > -\fep -1 = -K,
\end{equation*}
because $\fep^2 + \fep < \ep^2$ holds in this case.

\ref{item:prop_fbound}) For $\ell \in S_k$, item \ref{item:prop_fbound}) follows from item~\ref{item:prop_f}).
Otherwise, for $\ell \ge k - \underline{k} + 1$, obtain from the time derivative of \eqref{eq:worse_f_1} for $t \ge k - \underline{k} + 1$
\begin{equation}
    \tilde{\f}^{(1)}(t) 
    \ge \tilde{\f}^{(1)}(k - \underline{k} + 1) = \frac{4}{\ep^2} + \tilde{\f}^{(1)}(k - \underline{k}) \ge \frac{2}{\ep^2} > 0
\end{equation}
for $k > \underline{k}$ due to item~\ref{item:prop_dfboundk}) and
\begin{align}
    \tilde{\f}^{(1)}(t) &\ge \tilde{\f}^{(1)}(1) = -\frac{4 (k-1)}{\ep^2} + \bar{a}_k = -\frac{2 k - 4}{\ep^2} + \frac{2}{k} \nonumber\\
    &\ge -\frac{2 \ep - 2}{\ep^2} + \frac{2}{\ep + 1} = \frac{2}{\ep^2(\ep+1)^2} > 0
\end{align}
for $k = \underline{k} \le K \le \ep + 1$.
Hence, $\tilde{\f}$ is strictly increasing on $[k - \underline{k} + 1, k]$ and, since $\tilde{\f}(k) =1$, it suffices to show $\tilde{\f}(k-\underline k +1) \ge -1$.
To see this, assume the opposite
\begin{equation*}
    -1 >\tilde{\f}(k-\underline{k}+1) =-1+\frac{2}{\underline{k}}+\frac{2}{\ep^2}-\frac{2\underline{k}}{\ep^2}
\end{equation*}
  or equivalently that $\ep^2+\underline{k}-\underline{k}^2<0$. For $\underline{k} \le \fep$ this is impossible, because then $\ep \ge \underline{k}$; hence,  $\underline{k}=\fep+1=K$. Then,
  $\ep^2+\fep+1-(\fep+1)^2<0$ or equivalently
  $\ep^2<\fep^2+\fep$ which contradicts the fact that $K=\fep + 1$.

\ref{item:prop_dfbound}) Note that $\tilde{\f}^{(2)}(t)=4/\varepsilon^2$ for $t\in (k-\underline{k},k)$. The result is thus established if $\underline{k}= k$. Next, consider $k>\underline{k}= K$. From~\eqref{eq:worse_f_2}, $|\tilde{\f}^{(2)}(t)|=|2\tilde{\f}^{(1)}(\ell)|$ for $t\in(\ell-1,\ell)$ with $\ell\in S_k$. From item~\ref{item:prop_df}), then $|2\tilde{\f}^{(1)}(\ell)|=|2\tilde{\f}^{(1)}(k-\underline{k})|$ for $\ell\in S_k$, where $|2\tilde{\f}^{(1)}(k-\underline{k})|\leq 4/\varepsilon^2 = \tilde{L}$ due to item~\ref{item:prop_dfboundk}).
Thus, $|\tilde{\f}^{(2)}(t)|\leq \tilde{L}$ follows for almost every $t\in[0,k]$.
\end{IEEEproof}  
\subsection{The case with arbitrary positive $N,T,L$}
Let $\f(t):=N\tilde{\f}(t/T)$ for $t\in[0,kT]$ and  $\tilde{\f}(t)$ defined as in \eqref{eq:worse_f_1} and \eqref{eq:worse_f_2} with $\varepsilon=\frac{2}{T}\sqrt{\frac{N}{L}}$. First, Lemma~\ref{lem:prop_ftilde}, items~\ref{item:prop_f}) and~\ref{item:prop_df}), is used to conclude that $\f$ is continuously differentiable in $(0,kT)$. Next, Lemma~\ref{lem:prop_ftilde}\ref{item:prop_fbound}) is used to conclude that $\f(\ell T)\in[-N,N]$ for every integer $\ell \in [0,k]$. Moreover, Lemma~\ref{lem:prop_ftilde}\ref{item:prop_dfbound}) is used to conclude that for almost all $t\in[0,kT]$,
$$
\left| \f^{(2)}(t) \right| = \frac{N}{T^2} \left| \tilde{\f}^{(2)}(t/T) \right|\leq \frac{N}{T^2} \frac{4}{\ep^2} = \frac{N}{T^2} \frac{L T^2}{N} =L.
$$
Furthermore, using $\varepsilon=(2/T)\sqrt{N/L}$ in \eqref{eq:definitions} recovers the definitions in the statement of Theorem~ \ref{thm:zero-meas-worst} directly.

It follows that the function $\f$ satisfies \eqref{eq:fddotbnd}--\eqref{eq:basemnf} for some sequence $\{\eta_k\}$ and zero measurements, and hence $\X_k(\f) \in \C_k(\Oo)$ (recall Remark~\ref{rem:bgg}). In addition, $\f^1_k=\f^{(1)}(kT)=\bar h_k = N\overline{a}_k/T$. From~\eqref{eq:trueworst}, \eqref{eq:bounds} and~\eqref{eq:truelpineq}, then 
\begin{align}
    \label{eq:hklowerbound}
    \overline{f}_k^1(\Oo)\geq \overline\F_k^1(\Oo) \ge \f^1_k = \bar h_k.
\end{align}

Next, $\bar{h}_k$ is shown to be also an upper bound for $\overline{f}_k^1(\Oo)$.
\begin{lem}
\label{le:f_ak_bound}
Let $k,\underline k\in\N$ satisfy $k\ge \underline k$. Consider real numbers $f_j$, $f_j^1$ for $j=0,1,\ldots,k$ satisfying in \eqref{eq:cset} the inequalities corresponding to \eqref{eq:ff1const} and \eqref{eq:f1const} for $j=k-\underline k+1,\ldots,k$, and to \eqref{eq:fmbnd} for $j=k$ and $j=k-\underline k$. Let $\bar h_k$ be defined as in Theorem \ref{thm:zero-meas-worst}. Then, $f_k^1\leq \bar{h}_k$.
\end{lem}
\begin{IEEEproof}
From \eqref{eq:f1const} we know that $f_{j-1}^1\geq f_{j}^1-LT$. Using this relation repeatedly for $j=k,k-1,\ldots,k-\underline k+1$,
  \begin{align}
    \label{eq:topder1}
    f_{k-i}^1 &\ge f_k^1 - iLT\quad \text{for }i=1,\ldots,\underline k. 
  \end{align}
Similarly, from~\eqref{eq:ff1const}, we know that $f_{j-1}\leq f_j-f_{j}^1 T+T^2L/2$. Using this relation for $j=k,k-1,\ldots,k-\underline k + 1$ yields  
  \begin{align}
    \label{eq:topfk1}
    f_{k-i} &\le f_k - \sum_{j=0}^{i-1} f_{k-j}^1 T + iL\frac{T^2}{2},\quad i=1,\ldots,\underline k.
  \end{align}
Let $i=\underline k$ and use \eqref{eq:topder1} in \eqref{eq:topfk1} to obtain:
  \begin{align}
    f_{k-\underline k} &\le f_k - \sum_{j=0}^{\underline k-1} \left( f_k^1 - jLT \right) T + \underline kL \frac{T^2}{2}\notag\\
    \label{eq:f0ineq1}
    &\le f_k - \underline kT f_k^1 + LT^2 \underline k^2/2.
  \end{align}
  Using $-LT^2 \underline k^2/2 = 2N-\underline kT \bar h_k$ from the definition of $\bar h_k$,
  \begin{align}
  \label{eq:lower_fkk}
    f_k - f_{k-\underline k} \ge \underline kTf_k^1 - \frac{\underline k^2LT^2}{2} = 2N + \underline k T (f_k^1 - \bar h_k).
  \end{align}
  However, from \eqref{eq:fmbnd} we know that, $f_k\leq N$ and $-f_{k-\underline k}\leq N$. Thus, $f_k - f_{k-\underline k}\leq 2N$, which with \eqref{eq:lower_fkk} yields $f_k^1\leq \bar h_k$.
  \end{IEEEproof}

Combining \eqref{eq:hklowerbound} and Lemma~\ref{le:f_ak_bound} leads to $\overline{f}_k^1(\Oo) = \overline{\F}_k^1(\Oo) =\bar h_k$. From~\eqref{eq:cset}, it follows that $\x \in \C_k(\Oo) \Leftrightarrow -\x \in \C_k(\Oo)$. Therefore, it must happen that $\underline f_k^1(\Oo) = \underline{\F}_k^1(\Oo) = -\bar h_k$. Finally, recalling \eqref{eq:width} and \eqref{eq:truewidth}, then $w_k(\Oo) = \W_k(\Oo) = 2\bar h_k$. This establishes Theorem \ref{thm:zero-meas-worst}\ref{item:finconv}) and~\ref{item:0lp0fun}).

To prove item~\ref{item:noninc}), note that $\bar h_k$ in Theorem~\ref{thm:zero-meas-worst} satisfies $\bar h_k = (N/T) h(\underline k)$, with the latter defined as in \eqref{eq:definitions} and $\varepsilon^2 = 4N/LT^2$. Therefore, Theorem~\ref{thm:zero-meas-worst}\ref{item:noninc}) follows from Lemma~\ref{le:sequence_a}.

\subsection{The case with $\m_k\neq \Oo$}
The constraint set $\C_k(\m_k)$ has the following simple property, which will be instrumental in establishing Theorem~\ref{thm:zero-meas-worst}\ref{item:0worst}).
\begin{lem}
  \label{lem:Caffinv}
  Let $\m_k\in\R^{k+1}$ as in \eqref{eq:defmk} and $\x_k\in\R^{2k+2}$ with components named as in \eqref{eq:defxk} be such that $\x_k \in \C_k(\m_k)$. Let $\tilde\m_k\in\R^{k+1}$ have components $\tilde m_0, \tilde m_1,\ldots, \tilde m_k$ satisfying
  \begin{align}
    \label{eq:mtildeDef}
    \tilde m_j &= m_j + aj + b,\quad j=0,1,\ldots,k,
  \end{align}
  for some $a,b\in\R$ and define $\tilde\x_k := [ (\tilde f_{0:k})^T, (\tilde f_{0:k}^1)^T]^T$, with
  \begin{align*}
    \tilde f_j &:= f_j + aj + b, & \tilde f_j^1 &:= f_j^1 + a/T, \quad j=0,1,\ldots,k.
  \end{align*}
  Then, $\tilde\x_k \in \C_k(\tilde\m_k)$.
\end{lem}
\begin{IEEEproof}
  Directly from the definitions, it is clear that
  \begin{gather*}
    \tilde f_j - \tilde m_j = f_j - m_j,\quad
    \tilde f_j^1 - \tilde f_{j-1}^1 = f_j^1 - f_{j-1}^1,\\
    \tilde f_{j-1} - \tilde f_{j} + \tilde f_j^1 T = f_{j-1} - f_j + f_j^1 T.
  \end{gather*}
  Therefore, if \eqref{eq:ff1const}--\eqref{eq:fmbnd} are satisfied for $\x_k$ and $\m_k$, they will also be satisfied for $\tilde\x_k$ and $\tilde\m_k$. 
\end{IEEEproof}

Consider $\m_k$ with nonempty $\C_k(\m_k)$ and let $\tilde\m_k$ be defined as in \eqref{eq:mtildeDef}
with $a = (m_{k-\underline k} - m_k)/\underline{k}$ and $b = - m_{k} - a k$.
Then, $\tilde m_{k-\underline k} = \tilde m_k = 0$. By Lemma~\ref{lem:Caffinv}, it follows that
$\overline{f}_k^1(\tilde\m_k) = \overline f_k^1(\m_k) + a/T$ and $\underline f_k^1(\tilde\m_k) = \underline f_k^1(\m_k) + a/T$, so that $w_k(\tilde\m_k) = w_k(\m_k)$. 

Next, apply Lemma~\ref{le:f_ak_bound} to $\tilde \x_k \in \C_k(\tilde \m_k)$. This gives $\overline f_k^1(\tilde\m_k) \le \bar h_k$. By the symmetry of the constraints required by Lemma~\ref{le:f_ak_bound}, also $\underline f_k^1(\tilde\m_k) \ge -\bar h_k$. Therefore, 
\begin{align*}
  w_k(\Oo) = 2\bar h_k \ge w_k(\tilde\m_k) = w_k(\m_k) \ge \W_k(\m_k)
\end{align*}
for every $k$ and $\m_k$. Taking the supremum over all $\m_k$ yields
    $w_k(\Oo) \ge \bar w_k(L, N, T) \ge \bar\W_k(L,N,T) \ge \W_k(\Oo)$.
Theorem~\ref{thm:zero-meas-worst}\ref{item:0worst}) is then established recalling Theorem~\ref{thm:zero-meas-worst}\ref{item:0lp0fun}).

\section{Comparisons}
\label{sec:comparisons}

\newcommand{\Mlp}{M_{\mathrm{lp}}}
\newcommand{\Tlp}{T_{\mathrm{lp}}}
\newcommand{\Msm}{M_{\mathrm{sm}}}
\newcommand{\Tsm}{T_{\mathrm{sm}}}
\newcommand{\Mhg}{M_{\mathrm{hg}}}
\newcommand{\tauhg}{\tau_{\mathrm{hg}}}
\newcommand{\diffd}{\mathrm{d}}
\newcommand{\abs}[1]{\left|#1\right|}

This section compares the proposed differentiator's performance and accuracy to a linear high-gain and an exact sliding-mode differentiator.
For comparison purposes, each of those two differentiators is discretized using state-of-the-art techniques.
The proposed differentiator in Algorithm~\ref{algo:differentiator} is implemented by solving the linear programs using Yalmip~\cite{Lofberg2004YalmipMatlab} with the Matlab solver \texttt{linprog}.

Before doing the comparison, it is worthwile to note that Theorem~\ref{thm:zero-meas-worst} states the proposed differentiator's worst-case accuracy $\frac{1}{2} \bar{\W}_k(L, N, T)$ and the maximum time $K T$ it takes to achieve it.
For all values of $T$, $L$, and $N$, this accuracy is bounded from below by
\begin{equation}
    \frac{1}{2} \bar{\W}_k(L, N, T) \ge 2 \sqrt{NL}.
\end{equation}
This lower limit is also obtained for certain special combinations of $T$, $L$, $N$, as well as for $T \to 0$.
Exact differentiators have a similar inherent accuracy restriction, see \cite{Levant1998RobustTechnique,Levant2017SlidingApplication}.

\subsection{Linear High-Gain Differentiator}

In continuous time, a second order linear (high-gain) differentiator with identical eigenvalues and time constant $\tau$ is given by
\begin{align}
\label{eq:highgain:diff}
    \dot y_1 &= \frac{2}{\tau} (m - y_1) + y_2, &
    \dot y_2 &= \frac{1}{\tau^2} (m - y_1),
\end{align}
with output $\hat{\f}^1 = y_2$, input $m = \f + \eta$ and $\abs{\eta} \le N$.
From \cite{Vasiljevic2008ErrorObservers}, its optimal asymptotic accuracy is obtained as \mbox{$4 e^{-\frac{1}{2}} \sqrt{NL} \approx 2.43 \sqrt{NL}$}.
The corresponding optimal time constant is $\tau = e^{-\frac{1}{2}} \sqrt{N/L}$, which is hence chosen in the following.
For simulation purposes, the linear system \eqref{eq:highgain:diff} is discretized using the implicit Euler method.

\subsection{Robust Exact Sliding-Mode Differentiator}

As a sliding-mode differentiator, the robust exact differentiator proposed in \cite{Levant1998RobustTechnique} is used.
In continuous time, it is
\begin{subequations}
\label{eq:smdiff}
\begin{align}
      \dot y_1 &= k_1 \abs{m - y_1}^{\frac{1}{2}} \sign(m - y_1) + y_2, \\
      \dot y_2 &= k_2  \sign(m - y_1),
\end{align}
\end{subequations}
with output $\hat{\f}^1 = y_2$, input $m = \f + \eta$ with $\abs{\eta} \le N$ and positive parameters $k_1, k_2$.
It is discretized using the matching method proposed in \cite{Koch2018DiscreteDifferentiators} and simulated using the toolbox~\cite{Andritsch2021RobustFeatures}.
Parameters are selected as $k_1 = 2 r$ and $k_2 = r^2$, with robustness factor $r$ as in \cite{Andritsch2021RobustFeatures} set to $r = 1.5 \sqrt{L}$.

\subsection{Comparison}
\begin{figure}
    \centering
    \includegraphics{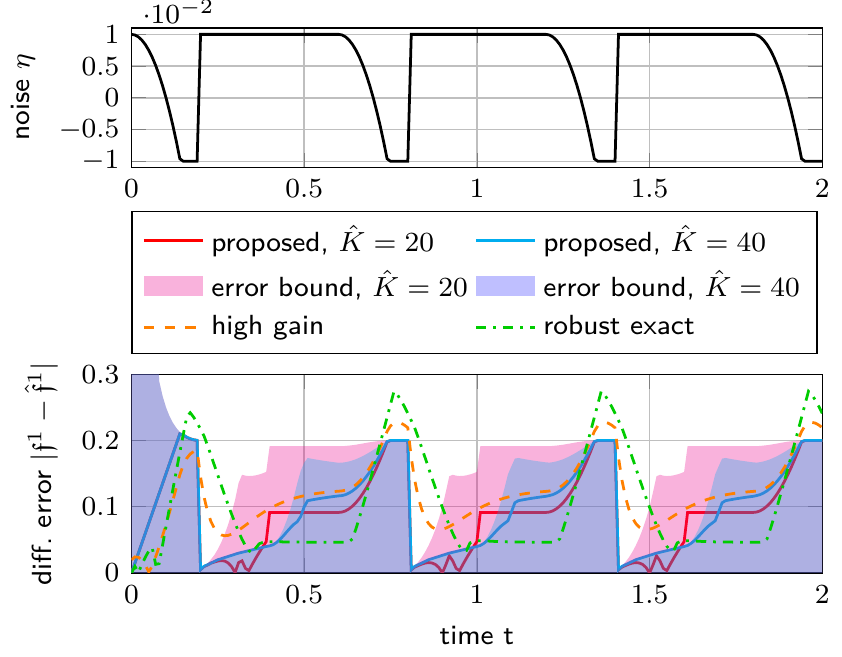}
    \caption{Bounded noise $\eta$ added to the signal $\f(t) = t^2/2$ sampled with \mbox{$T = 10^{-2}$}, and corresponding differentiation errors for proposed differentiator, linear high-gain differentiator, and  sliding-mode differentiator from a simulation with $N = 10^{-2}$ and $L = 1$. For the proposed differentiator, the error bounds obtained along with the estimate from the linear programs in Algorithm~\ref{algo:differentiator} are also shown.}
    \label{fig:simulation}
\end{figure}

For the comparison, the signal $\f(t) = L t^2/2$ and noise
\begin{equation*}
    \eta(t) = \begin{cases}
      \max(-N, N - L (t - c \lfloor \frac{t}{c} \rfloor)^2) & t - c \lfloor \frac{t}{c} \rfloor < 2 \sqrt{\frac{N}{L}} \\
      N & \text{otherwise}
    \end{cases}
\end{equation*}
with constant $c = 6 \sqrt{N/L}$ are sampled with $T = 10^{-2}$.
Parameters are selected as $L = 1$, $N = 10^{-2}$.
For these particular parameters, Theorem~\ref{thm:zero-meas-worst} yields $K = 20$ and an optimal worst-case accuracy $h_o(K) = 0.2$.

Fig.~\ref{fig:simulation} depicts the noise as well as the differentiation error of all differentiators.
For the proposed differentiator, two values of $\hat K$ are considered and the error bounds, i.e., the values of $(\overline{f}_k^1 - \underline{f}_k^1)/2$, as obtained from the linear program are shown as well.
One can see that, after an initial transient of duration $K T = 0.2$, the proposed differentiator achieves the best worst-case accuracy of $h_o(K) = 0.2$, as expected from the theoretical results.
Moreover, increasing $\hat K$ improves the error bound obtained along with the estimate.
The high-gain differentiator leads to a larger but smoother error overall.
The robust exact differentiator, finally, exhibits the largest worst-case errors, because it attempts to differentiate exactly also the noise, but is the most accurate one for constant noise.
\section{Conclusion}
\label{sec:conclusion}

A differentiator for sampled signals based on linear programming was proposed. It is shown that the best worst-case accuracy is obtained with a fixed number of discrete-time measurements, which allows limiting its computational complexity.
Comparisons to a linear high-gain differentiator and a standard sliding-mode differentiator exhibited a higher accuracy. However, depending on the sampling time, the increased accuracy comes at a higher computational cost. 

\appendix
To show that, as stated in Remark~\ref{rem:bgg}, a nonempty constraint set $\mathcal{C}_k$ does not necessarily imply existence of a function $\f$, consider $L = 2$,  $T = 1$, $N = 0$, $k = 2$ and measurements \mbox{$\m = [0, 0, 4]^T$}.
It is easy to check that every $\x \in \mathcal{C}_k(\m)$ has the form $\x = [0, 0, 4, f_0^1, 1, 3]^T$ with $f_0^1 \in [-1, 3]$.
By symmetry with respect to time reversal, any function $\f$ satisfying  \eqref{eq:fddotbnd} also has to satisfy $|\f_{j} - \f_{j+1} + \f_j^1 T| \le L \frac{T^2}{2}$ in addition to \eqref{eq:ff1const}.
Adding this inequality with $j =1$ to the constraints as $\abs{f_1 - f_2 + f^1_1} \le 1$ yields a contradiction.
Hence, no function $\f$ satisfying \eqref{eq:fddotbnd} exists for these measurements.


\begin{thebibliography}{10}

\bibitem{Reichhartinger2018SpecialDifferentiators}
M.~Reichhartinger, D.~Efimov, and L.~Fridman, ``{Special issue on
  differentiators},'' {\em Int. J. Control}, vol.~91, no.~9, pp.~1980--1982,
  2018.

\bibitem{RasoolMojallizadeh2021Discrete-timeAnalysis}
M.~R. Mojallizadeh, B.~Brogliato, and V.~Acary, ``{Discrete-time
  differentiators: design and comparative analysis}.'' {HAL} Preprint ID:
  hal-02960923, 2021.

\bibitem{Shtessel2014ObservationObservers}
Y.~Shtessel, C.~Edwards, L.~Fridman, and A.~Levant, ``{Observation and
  Identification via HOSM Observers},'' in {\em Sliding Mode Control and
  Observation}, pp.~251--290, Springer, 2014.

\bibitem{Rios2015fault}
H.~R{\'\i}os, J.~Davila, L.~Fridman, and C.~Edwards, ``Fault detection and
  isolation for nonlinear systems via high-order-sliding-mode
  multiple-observer,'' {\em Int. J. Robust Nonlin.}, vol.~25, no.~16,
  pp.~2871--2893, 2015.

\bibitem{Efimov2012ApplicationDetection}
D.~Efimov, L.~Fridman, T.~Ra{\"{i}}ssi, A.~Zolghadri, and R.~Seydou,
  ``{Application of interval observers and HOSM differentiators for fault
  detection},'' in {\em IFAC Proceedings}, vol.~8, pp.~516--521, 2012.

\bibitem{Bejarano2010HighInputs}
F.~J. Bejarano and L.~Fridman, ``{High order sliding mode observer for linear
  systems with unbounded unknown inputs},'' {\em Int. J. Control}, vol.~83,
  no.~9, pp.~1920--1929, 2010.

\bibitem{Vasiljevic2008ErrorObservers}
L.~K. Vasiljevic and H.~K. Khalil, ``{Error bounds in differentiation of noisy
  signals by high-gain observers},'' {\em Syst. Control Lett.}, vol.~57,
  no.~10, pp.~856--862, 2008.

\bibitem{Mboup2007AControl}
M.~Mboup, C.~Join, and M.~Fliess, ``{A revised look at numerical
  differentiation with an application to nonlinear feedback control},'' in {\em
  Mediterranean Conference on Control and Automation}, IEEE, 2007.

\bibitem{Othmane2020AnalysisDisturbances}
A.~Othmane, J.~Rudolph, and H.~Mounier, ``Analysis of the parameter estimate
  error when algebraic differentiators are used in the presence of
  disturbances,'' in {\em 21st IFAC World Congress}, pp.~572--577, 2020.

\bibitem{Levant1998RobustTechnique}
A.~Levant, ``{Robust Exact Differentiation via Sliding Mode Technique},'' {\em
  Automatica}, vol.~34, no.~3, pp.~379--384, 1998.

\bibitem{Levant2003}
A.~Levant, ``{Higher-order sliding modes, differentiation and output-feedback
  control},'' {\em Int. J. Control}, vol.~76, no.~9-10, pp.~924--941, 2003.

\bibitem{Fraguela2012}
L.~Fraguela, M.~T. Angulo, J.~A. Moreno, and L.~Fridman, ``{Design of a
  prescribed convergence time uniform Robust Exact Observer in the presence of
  measurement noise},'' in {\em Conference on Decision and Control},
  pp.~6615--6620, 2012.

\bibitem{Polyakov2019ConsistentSystems}
A.~Polyakov, D.~Efimov, and B.~Brogliato, ``{Consistent discretization of
  finite-time and fixed-time stable systems},'' {\em SIAM Journal on Control
  and Optimization}, vol.~57, no.~1, pp.~78--103, 2019.

\bibitem{Levant2013OnControl}
A.~Levant, ``{On fixed and finite time stability in sliding mode control},'' in
  {\em Conference on Decision and Control}, pp.~4260--4265, IEEE, 2013.

\bibitem{Carvajal-Rubio2021ImplicitDifferentiators}
J.~E. Carvajal-Rubio, J.~D. S{\'{a}}nchez-Torres, A.~G. Loukianov, M.~Defoort,
  and M.~Djemai, ``{Implicit and Explicit Discrete-Time Realizations of
  Homogeneous Differentiators},'' {\em Int. J. Robust Nonlin.}, 2021.

\bibitem{Koch2018DiscreteDifferentiators}
S.~Koch and M.~Reichhartinger, ``Discrete-time equivalent homogeneous
  differentiators,'' in {\em 15th International Workshop on Variable Structure
  Systems (VSS)}, pp.~354--359, IEEE, 2018.

\bibitem{Mazenc2011IntervalDisturbances}
F.~Mazenc and O.~Bernard, ``{Interval observers for linear time-invariant
  systems with disturbances},'' {\em Automatica}, vol.~47, no.~1, pp.~140--147,
  2011.

\bibitem{Guerra2017IntervalAccuracy}
M.~Guerra, C.~V{\'{a}}zquez, D.~Efimov, G.~Zheng, L.~Freidovich, and
  W.~Perruquetti, ``{Interval differentiators: On-line estimation of
  differentiation accuracy},'' in {\em European Control Conference},
  pp.~1347--1352, IEEE, 2017.

\bibitem{Lofberg2004YalmipMatlab}
J.~L\"ofberg, ``Yalmip: A toolbox for modeling and optimization in matlab,'' in
  {\em IEEE Int Conf Robot Autom}, pp.~284--289, 2004.

\bibitem{Levant2017SlidingApplication}
A.~Levant, M.~Livne, and X.~Yu, ``Sliding-mode-based differentiation and its
  application,'' in {\em 20th IFAC World Congress}, pp.~1699--1704, 2017.

\bibitem{Andritsch2021RobustFeatures}
B.~Andritsch, M.~Horn, S.~Koch, H.~Niederwieser, M.~Wetzlinger, and
  M.~Reichhartinger, ``The robust exact differentiator toolbox revisited:
  Filtering and discretization features,'' in {\em IEEE International
  Conference on Mechatronics (ICM)}, 2021.

\end{thebibliography}
\end{document}